# A GENERALIZATION OF OBRESHKOFF-EHRLICH METHOD FOR MULTIPLE ROOTS OF POLYNOMIAL EQUATIONS

## A. I. Iliev

**Introduction.** After 1960 the question of simultaneous finding all roots (SFAR) of polynomials became very actual and it is considered by many authors. The reason of this interest is the better behaviour of the methods for SFAR with respect to the methods for individual search of the roots. Methods for SFAR have a wider region of convergence and they are more stable. In several survey publications [1,2,3] this question is considered in details. The first methods for SFAR are related to the case when the roots are simple. The well-known method of Dochev [4] is for SFAR of algebraic polynomial with real and simple roots. The developments of this same method for the case of nonalgebraic polynomials (trigonometric, exponential and generalized) are performed in [5,6,7]. The classical method of Obreshkoff - Ehrlich [8] possessing cubic rate of convergence is also generalized [9]. Using the approach basing on the divided differences with multiple knots Semerdzhiev [10] generalized the method of Dochev to the case when the roots have arbitrary, but given multiplicities. The same question for the case of trigonometric and exponential polynomials is solved in [11,3]. The new methods preserve their quadratic rate of convergence. The method of Obreshkoff - Ehrlich is also generalized to the most general case [12,3] of polynomials upon some Chebyshev system, having multiple roots with given multiplicities. The rate of convergence is cubic but, unfortunately, this generalization requires at each iteration to calculate determinants, which is labour-consuming operation.

In this paper we develop a new method, which is a generalization of the Obreshkoff - Ehrlich method for the cases of algebraic, trigonometric and exponential polynomials. This method has a cubic rate of convergence. It is efficient from the computational point of view and can be used for SFAR if the roots have known multiplicities. This new method in spite of the arbitrariness of multiplicities is of the same complexity as the methods for SFAR of simple roots. We do not use divided differences with multiple knots and this fact does not lead to calculation of derivatives of the given polynomial of higher order, but only of first ones.

**Algebraic polynomials.** Let the algebraic polynomial

(1) $\quad f(x) = x^n + a_1 x^{n-1} + \ldots + a_n$

be given and $x_1, x_2, \ldots, x_m$ are his roots with given multiplicities $\alpha_1, \alpha_2, \ldots, \alpha_m$ respectively $(\alpha_1 + \alpha_2 + \ldots + \alpha_m = n)$. For SFAR of (1) we define the following iteration process

(2) $\quad x_i^{[k+1]} = x_i^{[k]} - \alpha_i f(x_i^{[k]}) \left[ f'(x_i^{[k]}) - f(x_i^{[k]}) Q_i'(x_i^{[k]}) / Q_i(x_i^{[k]}) \right]^{-1}$

$\quad i = \overline{1,m}, \quad k = 0, 1, 2, \ldots$

where
$$Q_i(x) = \prod_{j \neq i, j=1}^{m} \left(x - x_j^{[k]}\right)^{a_j}.$$

One can verify directly that the following lemma holds true

**Lemma 1.** If $a_1 = a_2 = ... = a_m = 1$ and $Q(x) = \prod_{j=1}^{m}\left(x - x_j^{[k]}\right)^{a_j}$ then
$$Q''\left(x_i^{[k]}\right) / Q'\left(x_i^{[k]}\right) = 2Q_i'\left(x_i^{[k]}\right) / Q_i\left(x_i^{[k]}\right).$$

So, in the case of simple roots the method (2) reduces into the Obreshkoff - Ehrlich method which can be written in the form

(3)
$$x_i^{[k+1]} = x_i^{[k]} - f\left(x_i^{[k]}\right)\left[f'\left(x_i^{[k]}\right) - \frac{1}{2}f\left(x_i^{[k]}\right)Q''\left(x_i^{[k]}\right) / Q'\left(x_i^{[k]}\right)\right]^{-1}$$

$$i = \overline{1,m}, \quad k = 0, 1, 2, ...$$

The cubic rate of convergence can be established by the following theorem

**Theorem 1.** Let $q$, $c$ and $d \stackrel{def}{=} \min_{i \neq j}|x_i - x_j|$ be real constants such that the following inequalities are satisfied

$1 > q > 0$, $c > 0$, $d - 2c > 0$, $0 < c^2(n - 3a_i) + c(n + (3d-1)a_i) < d^2 a_i$, $i = \overline{1,m}$.

If the initial approximations $x_1^{[0]},...,x_m^{[0]}$ to the exact roots $x_1,...,x_m$ of (1) are chosen so that the inequalities $\left|x_i^{[0]} - x_i\right| \leq cq$, $i = \overline{1,m}$ hold true then for every natural $k$ the inequalities $\left|x_i^{[k]} - x_i\right| \leq cq^{3^k}$, $i = \overline{1,m}$ also hold true.

The proof of Theorem 1 can be carried out by induction with respect to the number of the iteration $k$.

**Example 1.** For the equation $(x-2)^2(x-3)^3(x-5) = 0$ at the initial approximations $x_1^{[0]} = 0.4$, $x_2^{[0]} = 3.5$ and $x_3^{[0]} = 8$ using the formula (2) we receive the roots with 18 decimal digits after only 4 iterations.

**Trigonometric polynomials.** The presented form (3) of the Obreshkoff - Ehrlich method prompted the generalizations [9] of (3) to the case of trigonometric and exponential cases.

Namely, for the trigonometric polynomial

(4) $$T_n(x) = \frac{a_0}{2} + \sum_{l=1}^{n}(a_l \cos lx + b_l \sin lx)$$

where at least one of the leading coefficients $a_n$ and $b_n$ is not zero and which has real roots $x_1,...,x_m$ with multiplicities $a_1, a_2,...,a_m$ $\left(a_1 + a_2 + ... + a_m = 2n\right)$ we can use the iteration method

(5)
$$x_i^{[k+1]} = x_i^{[k]} - a_i T_n\left(x_i^{[k]}\right)\left[T_n'\left(x_i^{[k]}\right) - T_n\left(x_i^{[k]}\right)Q_i'\left(x_i^{[k]}\right) / Q_i\left(x_i^{[k]}\right)\right]^{-1}$$

$$i = \overline{1,m}, \quad k = 0, 1, 2, ...$$

where

$$Q_i(x) = \prod_{j \neq i, j=1}^{m} \sin^{a_j}\left((x - x_j^{[k]})/2\right).$$

The formula (5) at $a_1 = a_2 = \ldots = a_m = 1$ coincides with the analogue of Obreshkoff - Ehrlich formula [9] for trigonometric polynomials

**Theorem 2.** Let us denote $d \stackrel{def}{=} \min_{i \neq j} |x_i - x_j|$. Let $c$, $q$ and $x$ be positive real numbers so that $q < 1$, $2c < x$, $d - 2c > 0$ and $\max_{i \neq j} |x_i - x_j| < 2p - 2x$. Denote the expression $\min\left\{\left|\sin\frac{x}{2}\right|, \left|\sin\left(\frac{d}{2} - c\right)\right|\right\}$ by $A$. If $c^2\left(4n + a_i(9A^2/8 - 2)\right) < A^2 a_i$, $i = \overline{1,m}$ and initial approximations $x_i^{[0]}$, $i = \overline{1,m}$ are chosen so that $|x_i^{[0]} - x_i| \leq cq$, $i = \overline{1,m}$ then for every natural $k$ the inequalities $|x_i^{[k]} - x_i| \leq cq^{3^k}$, $i = \overline{1,m}$ also hold true.

**Example 2.** For the trigonometric polynomial
$$T_3(x) = \sin^3((x-1)/2) \sin^2((x-2)/2) \sin((x-2.5)/2)$$
at initial approximations $x_1^{[0]} = 0.2$, $x_2^{[0]} = 1.7$ and $x_3^{[0]} = 3$ we reach the roots of $T_3(x)$ with an accuracy of 18 digits at the $5^{th}$ iteration.

**Exponential polynomials.** Let us now consider the polynomial

(6) $$E_n(x) = \frac{a_0}{2} + \sum_{l=1}^{n}\left(a_l \, ch\, lx + b_l \, sh\, lx\right) = \frac{a_0}{2} + \sum_{l=1}^{n}\left(a'_l \, \ell^{lx} + b'_l \, \ell^{-lx}\right).$$

We suppose that at least one of the leading coefficients $a_n$ or $b_n$ is not zero and that $E_n(x)$ has real roots $x_1, x_2, \ldots, x_m$ with known multiplicities $a_1, a_2, \ldots, a_m$ ($a_1 + a_2 + \ldots + a_m = 2n$) correspondingly. The roots of (6) can be refined simultaneously with the help of the computational scheme

(7) $$x_i^{[k+1]} = x_i^{[k]} - a_i E_n(x_i^{[k]})\left[E'_n(x_i^{[k]}) - E_n(x_i^{[k]})Q'_i(x_i^{[k]})/Q_i(x_i^{[k]})\right]^{-1}$$
$$i = \overline{1,m}, \quad k = 0, 1, 2, \ldots$$

where

$$Q_i(x) = \prod_{j \neq i, j=1}^{m} sh^{a_j}\left((x - x_j^{[k]})/2\right).$$

**Theorem 3.** Denote $\min_{i \neq j}|x_i - x_j|$ by $d$. Let $q$ and $c$ be real numbers such that $1 > q > 0$, $c > 0$, $d - 2c > 0$, $c^2\left(4n + (S^2 - 2)a_i\right) < S^2 a_i$, $i = \overline{1,m}$, where $S$ replace the expression $sh((d-2c)/2)$. If the initial approximations $x_i^{[0]}$, $i = \overline{1,m}$ are taken such that $|x_i^{[0]} - x_i| \leq cq$, $i = \overline{1,m}$ then for every $k \in N$ the inequalities $|x_i^{[k]} - x_i| \leq cq^{3^k}$, $i = \overline{1,m}$ hold true.

**Example 3.** The iteration method (7) was applied for SFAR of the exponential polynomial $E_2(x) = sh^2((x+2)/2) sh^2((x-3)/2)$. Using initial approximations $x_1^{[0]} = -1$ and $x_2^{[0]} = 4$ by the formula (7) we receive the roots with 18 decimal digits after only 4 iterations.

*Acknowledgement*. The author would like to express his deep gratitude to Dr Khristo Semerdzhiev for formulating the problem and for his help during these investigations.

University of Plovdiv
http://www.pu.acad.bg
Faculty of Mathematics and Informatics
http://www.fmi.pu.acad.bg
24 Tzar Assen St.
Plovdiv 4000
e-mail: aii@pu.acad.bg
URL: http://anton.iliev.tripod.com